\documentclass[12pt,reqno]{amsart}
\usepackage{amsmath,amssymb,latexsym,textcomp,mathrsfs}
\usepackage[all]{xy}
\usepackage{graphicx}
\usepackage{bm,amsmath, amsthm, amssymb, amsfonts}
\usepackage{times}
\usepackage[english]{babel}

\setbox0=\hbox{$+$}
\newdimen\plusheight
\plusheight=\ht0
\def\+{\;\lower\plusheight\hbox{$+$}\;}

\setbox0=\hbox{$-$}
\newdimen\minusheight
\minusheight=\ht0
\def\-{\;\lower\minusheight\hbox{$-$}\;}

\setbox0=\hbox{$\cdots$}
\newdimen\cdotsheight
\cdotsheight=\plusheight
\def\cds{\lower\cdotsheight\hbox{$\cdots$}}

\numberwithin{equation}{section}
\theoremstyle{plain}

  \newenvironment{nouppercase}{%
   \renewcommand{\uppercasenonmath}[1]{}}{}
	
	 \newcommand{\Keywords}[1]{\par\noindent
   {\small{Keywords and phrases}: #1}}
   
   \newcommand{\AMS}[1]{\par\noindent
   {\small{AMS Subject Classification (2010)}: #1}}

\begin{document}

\title{ SEMI $\lambda^*$-CLOSED SETS AND  NEW SEPARATION AXIOMS IN ALEXANDROFF SPACES}

 \author{Amar Kumar Banerjee$^1$}
 \author{Jagannath Pal$^2$}
 \newcommand{\acr}{\newline\indent}
 \maketitle
 \address{{1\,} Department of Mathematics, The University of Burdwan, Golapbag, East Burdwan-713104,
 West Bengal, India.
 Email: akbanerjee1971@gmail.com\acr
 {2\,} Department of Mathematics, The University of Burdwan, Golapbag, East Burdwan-713104,
 West Bengal, India. 
 Email:jpalbu1950@gmail.com\\}   
\begin{abstract}
Here we have studied the idea of semi  $\lambda^*$-closed sets and investigate some of their properties in  spaces considered by A. D. Alexandroff \cite{AD}. We have introduced some new separation axioms  namely semi-$ T_\frac{\omega}{4} $,  semi-$ T_\frac{3\omega}{8} $, semi-$ T_\frac{5\omega}{8} $ and their inter-relation with semi-$ T_0 $ and semi-$ T_1 $.  We have shown that under certain conditions these axioms are equivalent.

\end{abstract}

\begin{nouppercase}
\maketitle
\end{nouppercase}

\let\thefootnote\relax\footnotetext{
\AMS{Primary 54A05, 54A10, 54D10}

}

\section{\bf Introduction}
\label{sec:int}
Topological spaces have been generalised in several ways. By weakening the union requirements A. D. Alexandroff (1940) \cite{AD} generalised the idea of topological space to a $ \sigma $-space (or simply space) by taking only countable unions of open sets to be open. Later many works on generalisation of sets were done in a more general structure of a bispace \cite{BS} and \cite{BP}. Also several topological properties were studied in a bispace in \cite{AS} and \cite{AM}. 

The idea of generalised closed sets in a topological space was given by Levine \cite{NL}. In 1987, Bhattacharyya and Lahiri \cite{BL} introduced the class of semi-generalised closed sets in a topological space. Many works on semi-generalised closed sets ($ sg $-closed sets) have been done \cite{DM}, \cite{MBD}, \cite{SMB},  etc. where more references can be found. By taking an equivalent form of $ sg $-closed sets,  P. Das and Rashid (2004) \cite{RD} gave the idea of a generalisation of semi-closed sets in the $ \sigma $-spaces called semi $ g^* $-closed sets and studied its various properties. Recently M. S. Sarsak (2011) \cite{MS} studied $ g_\mu $-closed sets in a generalised topological space   and introduced new separation axioms namely $ \mu$-$T_\frac{1}{4},  \mu$-$T_\frac{3}{8} $ and $ \mu$-$T_\frac{1}{2} $ axioms by defining $ \lambda_\mu $-closed sets and investigate their properties and relations among the axioms. In \cite{AJ} the idea of $ \lambda^* $-closed sets was given in a space.  

In this paper we have studied the idea of semi-generalised closed sets namely $ sg^*$-closed  sets in Alexandroff spaces. We also introduce the idea of semi $ \lambda^* $-closed sets in Alexandroff spaces and  investigate various properties of these sets and obtain new separation axioms like semi-$ T_\frac{\omega}{4}$, semi-$T_\frac{3\omega}{8}$ and semi-$ T_\frac{5\omega}{8} $ axioms in Alexandroff spaces. Some  examples are subtantiated where necessary to enrich it.

 \section{\bf Preliminaries}
 \label{sec:pre}

\textbf{Definition 2.1}  \cite{AD}: A non empty set $X$ is called a $ \sigma $-space or simply a space if in it is chosen a system of subsets $\mathcal{F}$ satisfying the following axioms:

(1)	 The intersection of a countable number of sets in $\mathcal{F}$ is a set in $\mathcal{F}$.     

(2)  The union of a finite number of sets in $\mathcal{F}$ is a set in $\mathcal{F}$

(3) The void set is a set in $\mathcal{F}$.

(4)	The whole set $ X $ is a set in $\mathcal{F}$.

Sets of $\mathcal{F}$ are called closed sets. Their complementary sets are called open sets. The collection of all such open sets will sometimes be denoted by $\tau $ and the space by    $(X ,\tau )$. When there is no confusion, the space  $(X , \tau)$ will simply be denoted by $ X $. The compliment of a set $ A $ is denoted by  $ A ^c $.

Note that a topological space is a space but in general $ \tau $ is not a topology as can be easily seen by taking $ X = R $ and $ \tau $ as the collection of all $F_\sigma$-sets in $ R $. Several examples of spaces are seen in \cite{PD}, \cite{DR}, \cite{LD}.The definition of closure of a set and interior of a set in a space are similar as in the case of a topological space. Note that closure of a set in a space may not be closed in general. Also interior of a set in a space may not be open. 

Throughout the paper $ X $ stands for a space and   sets are always subsets of $ X $ unless otherwise stated. The letters $ R $ and $ Q $ stand respectively for the set of real numbers and the set of rational numbers.\\

\textbf{Definition 2.2} \cite{LD} : A set $ A $ in $ X $ is said to be semi-open if there exists an open set $ G $ in $ X $ such that $ G\subset A \subset \overline{G} $.\\

\textbf{Definition 2.3} \cite{PR} : A set $ A $ in $ X $ is said to be semi-closed if and only if $ X - A $ is semi-open.\\

The class of all semi-open sets and semi-closed sets in $ X $ will be respectively denoted by s.o. $ (X) $ and s.c. $ (X) $.\\

 Note that obviously open set is semi-open and closed set is semi-closed.\\

\textbf{Definition 2.4} \cite{RD} :  Two non-void sets $ A, B $ in $ X $ are said to be semi-separated if there exist two semi-open sets $U,  V $ such that $A\subset U, B\subset V $ and $ A\cap V=B\cap U=\emptyset$.\\

\textbf{Definition 2.5} \cite{RD} : A  space $(X, \tau)$ is called  semi-$ T_0 $ if for any two distinct points $ x,y $ of $ X $, there exists a semi-open set $ U $ which contains one of the points not the other.

Observe that a space $(X,  \tau)$  is semi-$T_0$ if and only if for any pair of distinct points $x, y \in X $,  there is a set $A$ containing one of the points, say $x$, but not $y$ such that $ A $ is either semi-open or semi-closed.

\textbf{Definition 2.6} \cite{RD} : A space $(X, \tau)$ is said to be semi-$T_1$ space if for any two distinct points $x, y \in X$, there are semi-open sets  $ U, V $ such that $x \in U,  y \not \in U,  y\in V , x\not\in V$.\\

\textbf{Definition 2.7 \cite{PR}:} A point $ x\in X $ is said to be a semi-limit point of  $ A \subset X $ in a space $ (X, \tau) $ if and only if for any semi-open set $ U $ containing $ x $,  $ U \cap(A - \{x\})\not=\emptyset $ . The set of all semi-limit points of $ A $ is called the semi derived set of $ A $ and is denoted by $ A_s' $. \\

\textbf{Definition 2.8}  \cite {PR}: Semi closure of a set $ A $ is denoted by $ \overline{A}_s $ and is defined by $ \overline{A}_s = A\cup A_s' $ or equivalently $ \overline{A}_s = \cap \{F : F $ is  semi-closed containing $ A $ \}.\\

\textbf{Theorem 2.9 } \cite{PR}: Countable intersection of semi-closed sets in a space is semi-closed.\\

In \cite{RD} it is shown that countable union of semi-open sets is semi-open.\\

\section{\bf Semi $ g^* $-closed sets, Semi $ \wedge _\tau $-sets,  Semi generalised $ \wedge _\tau $-sets and Semi-$ T_\omega $ spaces.}

\textbf{Definition 3.1} \cite {BL} ): In a topological space, a subset $ A $ is said to be  semi-generalised closed set if and only if $ scl(A)\subset O $ whenever $ A \subset O $ and $ O $ is semi-open, where $ scl(A) $ denotes the semi-closure of $ A $.\\

\textbf{Definition 3.2 }  \cite {RD} ): A subset $ A $ of $ X $ is said to be  semi $ g^* $-closed set $(sg^*$-closed set in short) if and only if there is a 
semi-closed set $ F $ containing $ A $ such that $ F\subset O $ whenever $ A\subset O $ and $ O $ is semi-open.\\

\textbf{Definition 3.3}  (\cite{RD}): A set $A$ is called semi $ g^*$-open ($sg^*$-open)  if $ X - A $ is $s g^*$-closed.\\

\textbf{Remark 3.4 :} Clearly by definition, every semi-closed set is $ sg^*$-closed set, but the converse may not be true which has been shown by an Example 1 in \cite{RD}.\\

\textbf{Definition  3.5} (\cite{RD}) : For $ A\subset X $, the semi-kernel of $ A $ denoted by $ sA_\tau^\wedge $ is the set $ \cap \{U: U$ is semi-open, $U\supset A\} $.\\

\textbf{Definition  3.6} : If $ A\subset X $, then we define $ sA_\tau^\vee = \cup \{F: F$ is semi-closed, $F\subset A\} $.\\

\textbf{Theorem  3.7}  \cite{RD}: A set $ A $ is $ sg^* $-closed if and only if there is a semi-closed set $ F $ containing $ A $ such that $ F\subset sA_\tau^\wedge $. \\

\textbf{Definition  3.8 :} A set $ A $ in $ X $ is called a semi $ \wedge _\tau $-set ($ s\wedge _\tau $-set in short) if and only if $ A = sA_\tau^\wedge $ or equivalently, $ A $ is the intersection of all semi-open sets containing $ A $.\\

\textbf{Definition  3.9 :}  A set $ A $ in $ X $ is called a $ s\vee  _\tau $-set if and only if $ A = sA_\tau^\vee $ or equivalently, $X - A $ is $ s\wedge_\tau $-set i.e. $ A $ is the union of all semi-closed sets contained  in  $ A $.\\

It may be noted that semi-open sets may not be open, semi closure of a set may not be semi-closed, and finite intersection of semi-open sets may not be semi-open. For the sake of completeness
we are giving below some new examples to support the above assertions.\\

\textbf{Example 3.10} : Semi-open set may not be open.

Suppose $ X = R - Q $  and $ \tau = \{X, \emptyset, G_i\} $ where $ \{G_i\} $ are the countable subsets of $ X $ each containing $ \sqrt{3} $. Then $ (X, \tau) $ is a space but not a topological space. Then $ \overline{G_i} = X $ for each $ i $. Now let $ A $ be the set of all irrational numbers in $ (1, 2) $. Take $ B = \{\sqrt{2}, \sqrt{3}\} $, then $ B $ is an open set and $ \overline{B} = X $. Therefore $ B\subset A\subset \overline{B} $. So $ A $ is a semi-open set, but not an open set.\\

\textbf{Example 3.11 :} Semi closure of a set may not be semi-closed.

Suppose $ X = R - Q $ and $ \tau = \{X, \emptyset,  G_i\} $ where $\{ G_i\}$ are the  countable subsets of $ X $. Then $ (X, \tau) $ is a space but not a topological space. Let $ A $ be the set of all irrational numbers in $ (0, \infty) $. Then $ \overline{A}_s = A $. But $ A $ is not semi-closed, since $ X - A $, the set  of all irrational numbers in $ (-\infty, 0) $, is not semi-open.\\

\textbf{Example 3.12 :} Finite intersection of semi-open sets may not be semi-open.

Suppose $ X = R - Q $ and $ \tau = \{X, \emptyset, G_i\} $ where $\{G_i\} $  are countable subsets of $ X - \{\sqrt{2}\} $. Then $ (X, \tau) $ is a space but not a topological space. Let $ A = \{\sqrt{5}, \sqrt{3}, \sqrt{2}\}$. Then $ \{\sqrt{5}, \sqrt{3}\} \subset A \subset \overline{\{\sqrt{5}, \sqrt{3}\}} =  \{\sqrt{5}, \sqrt{3}, \sqrt{2}\}$ where $ \{\sqrt{5}, \sqrt{3}\}$ is an open set. This implies $ A $ is semi-open. Similarly the set $ B = \{\sqrt{11}, \sqrt{7}, \sqrt{2}\} $ is semi-open. But the set $ A\cap B = \{\sqrt{2}\} $ is not semi-open. \\

\textbf{Theorem  3.13}   \cite{RD}  : A subset $ A $ in $ X $ is $ sg^* $-closed if and only if there is a semi-closed set $ F $ containing $ A $ such that $ F - A $ does not contain any non-void semi-closed set.\\

\textbf{Theorem  3.14} (cf. Theorem 10 \cite{RD}): For each $ x \in X $, $ \{x\} $ is semi-closed or $ X - \{x\} $ is $ sg^* $-closed.\\

\textbf{Definition 3.15 :}    A set $A$ is called a semi generalised $\wedge_\tau$-set  ($sg\wedge_\tau$-set in short) if $sA_\tau^\wedge\subset F$ whenever $F\supset A$ and $F$ is semi-closed. 

A set $A$ is called a semi generalised $ \vee_\tau$-set  ($s g\vee_\tau $-set in short) if $ X - A  $ is $s g\wedge_\tau$-set.\\

\textbf{Theorem  3.16 :} For each $ x \in X $, $ \{x\} $ is either semi-open or $ sg\vee_\tau $-set.
\begin{proof}
Suppose $ x\in X $ and $ \{x\} $  is not semi-open, then $ X-\{x\} $ is not semi-closed. Then $ s(X - \{x\})_\tau^\wedge=\cap \{U, U $ is semi-open, $ U\supset (X -\{x\})\} \subset X $ where $ X $ is the only semi-closed set containing $ X - \{x\} $. Therefore  $ X-\{x\} $ is $s g\wedge_\tau $-set and so \{x\} is $s g\vee_\tau $-set.
\end{proof}

\textbf{Lemma  3.17} : Suppose $ A,  B $ be subsets of $ X $, the following can be easily verified.\\

(1) $ s\emptyset _\tau^\wedge = \emptyset, \quad  s\emptyset_\tau^\vee = \emptyset, \quad sX_\tau^\wedge = X, \quad sX_\tau^\vee = X $. \\

(2) $ A \subset sA_\tau^\wedge $, \quad  $ A\supset sA_\tau^\vee$\\

(3)  $ s(sA_\tau^\wedge)_\tau^\wedge = sA_\tau^\wedge $, \quad  $s(sA_\tau^\vee)_\tau^\vee = sA_\tau^\vee $\\

(4)  If $ A \subset B $ then $ sA_\tau^\wedge \subset sB_\tau^\wedge $.\\

(5)  If $ A \subset B $ then $ sA_\tau^\vee \subset sB_\tau^\vee $.\\

(6)  $ s(X\setminus A)_\tau^\wedge = X\setminus sA_\tau^\vee $, \quad  $ s(X\setminus A)_\tau^\vee = X\setminus sA_\tau^\wedge $.\\

\textbf{Theorem 3.18 :} If a set $ A $ is $ s\wedge_\tau $-set, then $ A $ is $ sg^* $-closed if and only if $ A $ is semi-closed.
\begin{proof}
Let $ A $ be $ s\wedge_\tau $-set and $ sg^* $-closed set. Then by Theorem 3.7, there exists a semi-closed set $ F $ containing $ A $ such that $ F\subset sA_\tau^\wedge =A $ i.e. $ F\subset A $. So $ F = A $ and so $ A $ is semi-closed. On the other hand, a semi-closed set is obviously a $ sg^* $-closed. Hence the result. 
\end{proof}

In particular, $ s(sA_\tau^\wedge)_\tau ^\wedge =  sA_\tau^\wedge$ and so $ sA_\tau^\wedge $ is a $ s\wedge_\tau $-set. Therefore $ sA_\tau^\wedge $ is $s g^* $-closed if and only if $ sA_\tau^\wedge $ is semi-closed.\\

\textbf{Theorem 3.19 :} If $ A\subset X $ and $ sA_\tau^\wedge $ is $ sg^* $-closed, then $ A $ is $ sg^* $-closed.
\begin{proof}
Suppose $ A \subset X $ and  $ sA_\tau^\wedge $ is $ sg^* $-closed, then there exists a semi-closed set $ F $ containing $ sA_\tau^\wedge $  such that $ F\subset s(sA_\tau^\wedge)_\tau^\wedge = sA_\tau^\wedge $. Since  $ F\supset sA_\tau^\wedge \supset A, A $ is $ sg^* $-closed.
\end{proof}

But the converse of the above theorem may  not be true as shown in the following example.\\

\textbf{Example 3.20 :} Suppose that $X = R - Q$  and $\tau = \{X, \emptyset,  G_i, A_i\} $ where $ \{G_i\} $ are the countable subsets of $ X $ containing $ \sqrt{2}$, $\{ A_i\} $ are the cocountable subsets of  $ X $ also containing $ \sqrt{2} $. Then $(X, \tau)$ is a space but not a topological space.  Let $A$ be a countably infinite subset of $X$ excluding $\sqrt{2}$.  Then $A$ is a closed set, so a semi-closed set which implies that $A$ is $sg^*$-closed. Again $ A $ is not an open set and also not a semi-open set since there is  no open set contained in $ A $. Now $sA_\tau^\wedge = \{ \sqrt{2}\}\cup A $ which is  an open set so a semi-open set. But $sA_\tau^\wedge $ is not a semi-closed set, since $ \sqrt{2}\in  sA_\tau^\wedge $.  Since $s A_\tau^\wedge $ is a $s \wedge_\tau $-set,  $sA_\tau^\wedge$ is not $sg^*$-closed, by Theorem 3.18.

\textbf{Theorem 3.21 :} Arbitrary union of $ s\vee _\tau $-sets is a $ s \vee_\tau $-set:

\begin{proof}
 Suppose $ A_i $ are $ s \vee_\tau $-sets, $ i\in I $ where $ i $ is an index set and $ A=\cup\{A_i: i\in I \} $. So $ A_i\subset A $ for each i. Therefore   $ sAi_\tau^\vee \subset sA_\tau^\vee $ for all $ i\in I $ and so $ \cup \{sAi_\tau^\vee: i\in I\}\subset sA_\tau^\vee $. So, $ A = \cup \{Ai: i \in I\}= \cup\{sAi_\tau^\vee: i\in I\}\subset sA_\tau^\vee  \subset A $ by Lemma 3.17 (2). Therefore $ A = sA_\tau^\vee $.
\end{proof}

\textbf{Theorem 3.22}: Intersection of two  $ s\vee _\tau $-sets is a $ s\vee_\tau $-set:
\begin{proof}
 Let $A , B$ be two $ s\vee_\tau $-sets , then $ A = sA_\tau^\vee, B = sB_\tau^\vee $.  Now $ A\cap B \subset A $ and also $ A\cap B \subset B $. So $ s(A\cap B)_\tau^\vee \subset sA_\tau^\vee $ and $ s(A\cap B)_\tau^\vee \subset sB_\tau^\vee $. Therefore  $ s(A\cap B)_\tau^\vee \subset  s A_\tau^\vee\cap sB_\tau^\vee $....(1).
 
  On the other hand, suppose that $ x\in sA_\tau^\vee\cap sB_\tau^\vee $ then $ x\in sA_\tau^\vee $ and $x \in sB_\tau^\vee$.  So there exist semi-closed sets $F , P$ such that $ x\in F\subset A, x\in P \subset B $. Therefore  $ x\in F\cap P \subset A\cap B $. This implies that $ x\in s(A\cap B)_\tau^\vee$ since $ F\cap P $ is a semi-closed set.  Therefore $ sA_\tau^\vee \cap sB_\tau^\vee \subset s(A\cap B)_\tau^\vee $.......(2). So by (1) and (2) $   s(A\cap B)_\tau^\vee =  sA_\tau^\vee \cap sB_\tau^\vee = A\cap B $.
\end{proof}

\textbf{Corollary 3.23} :  Let $(\ X, \tau)$ be a space. Then the collection of all $ s\vee_\tau $-sets forms a topology.

The proof follows from above two Theorems 3.21 and 3.22 and Lemma 3.17.\\

\textbf{Lemma 3.24 }: If $A , B$ are two subsets of $X$, then $s (A\cup B)_\tau^\wedge= sA_\tau^\wedge\cup sB_\tau^\wedge$.

\begin{proof}
Let $A , B$ be two subsets of $X$. Then  $ A\subset A\cup B $ implies that $s A_\tau^\wedge\subset s(A\cup B)_\tau^\wedge $ and $ B\subset A\cup B $ implies that $s B_\tau^\wedge\subset s(A\cup B)_\tau^\wedge $. Therefore ($s A_\tau^\wedge\cup sB_\tau^\wedge)\subset s(A\cup B)_\tau^\wedge $. Again, $s A_\tau^\wedge=\cap \{U_i: U_i\supset A,  U_i $ is semi-open\} and $s B_\tau^\wedge=\cap \{V_j:  V_j\supset B,  V_j$ is semi-open\}.  Therefore $s A_\tau^\wedge\cup sB_\tau^\wedge = \cap\{(U_i\cup V_j): U_i \supset A, V_j\supset B; U_i, V_j $ are semi-open\} $\supset\cap \{G: G\supset A\cup B, G$ is semi-open \}=$s (A\cup B)_\tau^\wedge $. Therefore $s A_\tau^\wedge\cup sB_\tau^\wedge=s(A\cup B)_\tau^\wedge $.
\end{proof}

But arbitrary union of $s \wedge_\tau $-sets is not always a $s \wedge_\tau $-set as can be seen from the example given below.\\

\textbf{Example 3.25} :Let $ X = R-Q $ and $ \tau $ = $\{X,  \emptyset, G_i\}$ where $\{ G_i\} $ be the countable subsets of $X$. Then $ (X, \tau) $ is a space but not a topological space. Here every singleton is an open set, so a semi-open set and hence a $s \wedge_\tau $-set. Now the set $ A=[1, 2]- Q $ is not an open set and also not a semi-open set because the closure of any open set $ B $ contained in $ A $ is $ B $. Again any set $ C (\not=X) $ containing $ A $ can not be a semi-open set, since closure of any open set $ D $ contained in $ C $ is $ D $. Therefore  
$s A_\tau^\wedge=X\not=A $. So $ A $ is not a $s \wedge_\tau $-set. But $ A=\cup \{\{r\}: r\in A\} $ where $ \{r\} $ is $s \wedge_\tau $-set.\\

\textbf{Theorem  3.26:} Intersection of two $s \wedge_\tau $-sets is a $s \wedge_\tau $-set.

\begin{proof}
Let $A , B$ be two $s \wedge_\tau $-sets , then $ A=sA_\tau^\wedge, B=sB_\tau^\wedge $.  Now $ A\cap B \subset A$ implies $ s(A\cap B)_\tau^\wedge\subset sA_\tau^\wedge $ and $ A\cap B \subset B$ implies $s(A\cap B)_\tau^\wedge\subset sB_\tau^\wedge $. So $s  (A\cap B)_\tau^\wedge \subset sA_\tau^\wedge\cap sB_\tau^\wedge=A\cap B \subset s(A\cap B)_\tau^\wedge$ by Lemma 3.17.  Therefore $  A\cap B = s(A\cap B)_\tau^\wedge $.The result follows.  
\end{proof}

\textbf{Definition 3.27}  \cite{RD}:  A space $ (X, \tau) $ is said to be a  semi-$T_\omega $ space  if and only if every $ s g^*$-closed set is semi-closed.\\

\textbf{Definition  3.28}(  \cite{RD}) : A topological spsce is called a semi-$ T_\frac{1}{2} $ space if and only if every $sg$-closed set is semi-closed.\\

 It is seen in \cite{RD} that semi-$ T_\frac{1}{2} $ axiom lies between semi-$ T_0 $ and semi-$ T_1 $ axioms. But semi-$ T_\omega $ axiom in a space does not have this property as evident from the Example 5 \cite{RD}. 

In Theorem  15  \cite{RD}, it is shown that every semi-$ T_\omega $ space is semi-$ T_0 $ space. But the converse is not true as shown in the Example 5 \cite {RD}. Also, in Examples 5 and 6 \cite{RD} it has been shown  that in general, semi-$ T_w$ and semi-$ T_1 $ axioms in a space are independent of each other.\\
 
\textbf{Theorem 3.29}  \cite{RD}): In a space, semi-$ T_1 $ axiom implies semi-$ T_\omega $ if the condition $ (P) $: ``Arbitrary intersection of closed sets is semi-closed in $ X $" is satisfied.\\

Next we see some necessary sufficient conditions for a space to be semi-$ T_\omega $. For this we introduce the following definition.\\

\textbf{Definition 3.30}  \cite{RD}) : For any $ E \subset X $,  let $ \overline{E_s^*}=\cap\{A: E\subset A ,$ A is $s g^*$-closed  set in X\},  then $\overline{E_s^*}$ is called $ sg^*$-closure  of the set $E$.\\

We denote the following sets as $ \mathcal {B} $  and $ \mathcal {B'} $ which will be used frequently in the sequel.

(i) $ \mathcal {B} = \{A: \overline{(X -A)_s}$ is semi-closed\}  and

(ii) $ \mathcal {B'} = \{A: \overline{(X -A)_s^*}$ is $sg^*$-closed\}. \\
 
\textbf{Theorem 3.31}  \cite{RD}  ):  A space $ (X,  \tau) $ is semi-$ T_\omega $ if and only if 

(a)  for each $ x \in X $,  $ \{x\} $ is either semi-open or semi-closed and

(b)  $ \mathcal {B} =  \mathcal {B'}$  

where $ \mathcal {B} = \{A: \overline{(X -A)_s}$ is semi-closed\}  and $ \mathcal {B'} = \{A: \overline{(X -A)_s^*}$ is $sg^*$-closed\}.

\section{\bf Semi $ \lambda^*$-closed sets and Semi $ \lambda^* $-open sets in  space:}

\textbf{Definition 4.1 :}  A subset $A$ of a space $ (X,  \tau) $  is said to be semi $ \lambda^*$-closed ($ s\lambda^* $-closed in short) if $ A = K\cap \overline{P}_s $  where $ K $ is a $s \wedge_\tau $-set   and $ P $ is a subset of $ X $.\\

\textbf{Definition 4.2 :} A subset $ A $ of a space $ (X,  \tau) $  is said to be semi $ \lambda^*$-open ($ s\lambda^* $-open in short) if $ X -A $ is  $ s\lambda^* $-closed.\\

\textbf{Definition 4.3:} The semi-interior of a set $ A $ in $ X $ is defined as the union of all semi-open sets contained in $ A $ and is denoted by $ sInt(A) $.\\

\textbf{Theorem 4.4 :} Suppose $ A\subset X $, then $ X - \overline{(X -A)}_s = sInt(A) $.
\begin{proof}
Suppose $ G $ be any semi-open set and $ G\subset A $. Then $ X - G \supset X -A $ implies $ \overline{(X-G)}_s \supset \overline{(X-A)}_s $ which implies that $ X-G \supset \overline{(X-A)}_s $. Therefore, $ G\subset  X - \overline{(X -A)}_s $ which implies that $ sInt(A) \subset X - \overline{(X -A)}_s $.

Conversely, $ \overline{(X -A)}_s = \cap\{P;  P $ is semi-closed, $ P\supset (X-A)\}$. So $ X - \overline{(X -A)}_s = \cup\{V; V $ is semi-open, $ V \subset A\}$. Each semi-open set $ V $ of the collection $ \{V = P^c: V $ is semi-open,$ V\subset A $\} is contained in $ sInt A $. Hence $ X - \overline{(X -A)}_s \subset sInt(A) $. Therefore, $ X - \overline{(X -A)}_s = sInt(A) $.
\end{proof}

\textbf{Lemma 4.5 :} Let   $ A $ be a subset of a space $ (X,  \tau) $, then $ \overline{\overline{A}}_s = \overline{A}_s $.
 \begin{proof}
 Since $ A \subset \overline{A}_s $ then $ \overline{A}_s\subset \overline{\overline{A}}_s$. Now $ \overline{\overline{A}}_s = \cap_{i\in\Lambda}\{ V_i; V_i $'s are semi-closed  and $ V_i \supset\overline{A}_s\} $. Again, $ \overline{A}_s$ = $\cap _{j\in\Lambda'}\{ U_j; U_j $'s are semi-closed  and $ U_j \supset A \} $. The sets $ U_j, j \in \Lambda' $ also contain $ \overline{A}_s $ i.e. $ U_j\supset \overline{A}_s , j \in \Lambda'$. Therefore the collection $ \{U_j, j\in\Lambda'\} \subset \{V_i, i\in\Lambda\} $ and so $ \cap _{j\in\Lambda'} U_j \supset\cap_{i\in\Lambda} V_i $ implies $ \overline{A}_s \supset \overline{\overline{A}}_s $. Hence $ \overline{\overline{A}}_s = \overline{A}_s $.
  \end{proof}

\textbf{Theorem 4.6 :} For a subset $ A $ of a space $ (X,  \tau) $ the following are equivalent:

(i) $ A $ is $ s\lambda^* $-closed

(ii)  $ A = sA_\tau^\wedge\cap \overline{P}_s, P \subset X $.

(iii)  $ A = K \cap \overline{A}_s, K $ is a  $ s\wedge_\tau $-set.

(iv)  $ A =   sA_\tau^\wedge\cap \overline{A}_s $.
\begin{proof}
(i) $\Longleftrightarrow $ (ii) : Suppose $ A $ is $ s\lambda^* $-closed, then $ A = K\cap \overline{P}_s $  where $ K $ is a $s \wedge_\tau $-set   and $ P \subset X $. Since $ A \subset K $ and $ A \subset \overline{P}_s, sK_\tau^\wedge \supset sA_\tau^\wedge $. Therefore $ A \subset sA_\tau^\wedge \cap \overline{P}_s \subset sK_\tau^\wedge \cap \overline{P}_s = K \cap \overline{P}_s = A $. So $ A =   sA_\tau^\wedge\cap \overline{P}_s $. Conversely, let $ A = sA_\tau^\wedge\cap \overline{P}_s, P \subset X $. Since $ sA_\tau^\wedge $ is a $ s\wedge_\tau $-set, then $ A $ is $  s\lambda^* $-closed. 

(i) $\Longleftrightarrow$ (iii) : Suppose $ A $ is $ s\lambda^* $-closed, then $ A = K\cap \overline{P}_s $  where $ K $ is a $ s\wedge_\tau $-set   and $ P \subset X $. Since $ A \subset K $ and $ A \subset \overline{P}_s $ implies $ \overline{A}_s \subset \overline{\overline{P}}_s = \overline{P}_s$. Therefore $ A \subset K \cap \overline{A}_s \subset K \cap \overline{P}_s = A $. So $ A = K\cap\overline{A}_s $. Conversely, suppose $ A = K\cap \overline{A}_s $  where $ K $ is a $ s\wedge_\tau $-set   and $ A \subset X $. Hence $ A $ is $  s\lambda^* $-closed. 

(i) $\Longleftrightarrow $ (iv) : Suppose $ A $ is $ s\lambda^* $-closed, then $ A = K\cap \overline{P}_s $  where $ K $ is a $ s\wedge_\tau $-set   and $ P \subset X $. Since $ K\supset A $ and $ A\subset \overline{P}_s $ then $ sK_\tau^\wedge \supset sA_\tau^\wedge $ and $ \overline{A}_s \subset \overline{\overline{P}}_s = \overline{P}_s$. Therefore $ A \subset sA_\tau^\wedge \cap \overline{A}_s \subset sK_\tau^\wedge \cap \overline{P}_s = K \cap \overline{P}_s = A $. So $ A = sA_\tau^\wedge\cap\overline{A}_s $. Conversely, let $ A = sA_\tau^\wedge\cap\overline{A}_s $ where $ A \subset X $ and $ sA_\tau^\wedge $ is a $ s\wedge_\tau $-set. Then $ A $ is $  s\lambda^* $-closed.
\end{proof}

\textbf{Remark 4.7 :} From Theorem 4.6 (iv)  we can say that a subset $A$ is said to be $s \lambda^*$closed   if $A$ can be expressed as the intersection of all semi-open sets and all semi-closed sets containing it.\\

Now we show by citing Examples 4.8 and 4.9 given below that there does not exist any relation between $ sg^* $-closed sets and $ s\lambda^* $-closed sets.\\

\textbf{Example 4.8 :} Example of a $ s\lambda^* $-closed set which is not $ sg^* $-closed.

Suppose $ X=R-Q $  and $ \tau = \{X, \emptyset, G_i\} $ where $ \{G_i\} $ are the countable subsets of $ X $, each contains $ \sqrt{3} $ then $ (X, \tau) $ is a space but not a topological space. Let $ A= \{\sqrt{3}, \sqrt{7}\} $. Therefore $ A $ is an open set so a semi-open set. Then $ sA_\tau^\wedge =A $ and $ \overline{A}_s = X $. Therefore $ A $ is  $ s\lambda^* $-closed. Here $ X $ is the only semi-closed set  containing  $ A $ which implies that $ A $ is not $ sg^* $-closed.\\

\textbf{Example 4.9 :} Example of a set $ A $ is $sg^* $-closed set which is not $ s\lambda^* $-closed.

Suppose $ X=R-Q, \{G_i\} $ are the countable subsets of $ X - \{\sqrt{3}\} $ and $ \tau = \{X, \emptyset, G_i\} $, then $ (X,  \tau) $ is a space but not a topological space. Let the subset $ A = X - \{\sqrt{3}\}$. Since $ \{\sqrt{3}\} $ is not a semi-open set, $ X - \{\sqrt{3}\} $ is not semi-closed. So $ \overline{A}_s = X $. Since closure of any open set $ B \subset A $ is $ B\cup \{\sqrt{3}\} $, there is no semi-open set except $ X $  containing $ A $. So $ sA_\tau^\wedge = X $, hence, $ sA_\tau^\wedge \cap \overline{A}_s = X \not= A $. Therefore $ A $ is not $ s\lambda^* $-closed. Again, since $ X $ is the only semi-open set containing $ A, A $ is  $sg^* $-closed.\\

\textbf{Example 4.10 :} Example of a  finite set which is not $ s\lambda^* $-closed.

Suppose $ X=R-Q $ and $ \tau = \{X, \emptyset, G_i\cup \{\sqrt{3},\sqrt{5}\} $ where $ \{G_i\} $ are the countable subsets of $ X $ containing $ \{\sqrt{3},\sqrt{5}\} $, then $ (X, \tau) $ is a space but not a topological space. Let $ A= \{\sqrt{3}, \sqrt{7}\} $. So $ A $ is not a semi-open set because no proper open set is contained in $ A $. But $ sA_\tau^\wedge =A\cup \{\sqrt{5}\}$ and $ \overline{A}_s=X $. Therefore  $ sA_\tau^\wedge\cap \overline{A}_s = (A \cup\{\sqrt{5}\})\cap X = A\cup \{\sqrt{5}\} \not= A $. Hence $ A $ is not a $ s\lambda^* $-closed set. \\

Clearly every $s \wedge _\tau $-set  is $s \lambda^*$-closed and semi-closed set is $s \lambda^*$-closed. But the converse may not be true as shown in the following example.\\

\textbf{Example 4.11 :} Let $X=R-Q$,  $ \tau = \{X,  \emptyset,  G_i\}$,  where $\{ G_i\} $  are the all countable subsets of $X$. So $ (X,  \tau) $ is a space but not a topological space.  Let $A$ be the set of all irrational numbers in $ (0,  \infty) $. Since $ sA_\tau^\wedge = X $ and $ \overline{A}_s = A $, $ A $ is $ s\lambda^* $-closed.  But $ A $ is not semi-closed because $ X - A $ is not semi-open. Also $ A $ is not a $ s\wedge_\tau $-set since  $ sA_\tau^\wedge = X \not = A $.\\

\textbf{Theorem 4.12:} A space $ (X,  \tau) $ is semi-$ T_\omega $ if and only if  every subset of $(X,\tau)$ is $s\lambda^*$- closed and $ \mathcal {B} =  \mathcal {B'} $ where $ \mathcal {B} $ and $  \mathcal {B'} $ are collection of sets as in Theorem 3.31.

\begin{proof} 
 Suppose every subset of  $(X,\tau)$ is $s\lambda^*$-closed and $ \mathcal {B} =  \mathcal {B'} $. Let $ x \in X $. We  claim that $\{x\}$ is either semi-open or semi-closed. Suppose $\{x\}$ is not semi-open, then $X-\{x\}$ is not semi-closed. Since $X-\{x\}$ is also a $s\lambda^*$-closed set then $X-\{x\} = s(X-\{x\})_\tau^\wedge\cap(\overline{X-\{x\}})_s$ [by Theorem 4.6 (iv)] = $s(X-\{x\})_\tau^\wedge\cap X = s(X- \{x\})_\tau^\wedge$. Therefore $X-\{x\}$ is a $s\wedge_\tau $-set. So $X-\{x\}$ is a semi-open set  which implies $\{x\}$ is semi-closed. Then by Theorem 3.31, $(X,\tau)$ is semi-$T_\omega$ space.
 
 Conversely, suppose that $ (X, \tau) $ is semi-$T_\omega$ space and $ A\subset X $.  Then by Theorem 3.31, every singleton is either semi-open or semi-closed and  $ \mathcal {B} =  \mathcal {B'} $.  So for each $ x\in X - A $,  either $ \{x\}\in s.o.(X) $  or  $ (X - \{x\}) \in s.o.(X) $. Let $ A_1=\{x : x\in X - A, \{x\}\in s.o.(X) \} $,  $ A_2=\{x : x\in X - A, X - \{x\}\in s.o.(X) \}$,  $ K=\cap[X - \{x\} : x\in A_2]=X - A_2 $ and   $P=\cap[X - \{x\} : x\in A_1]=X - A_1 $. 
 
 Note that $ K $ is a $s \wedge_\tau $-set   i.e. $ K=sK_\tau^\wedge$  and $ P=\overline{P}_s $.
 
 Now,  $K\cap\overline{P}_s=(X - A_2)\cap (X - A_1)=X - (A_1\cup A_2)=X -(X - A)=A$.  Thus $A$ is  $s\lambda^*$-closed. 
\end{proof}

\textbf{Remark 4.13 :}  It is already seen that if a subset $A$  is semi-closed then $A$ is $s g^*$-closed  and $ s \lambda^*$-closed . 

 But the converses of this result may not hold in a space as shown in the Example 1 \cite{RD} and in the Example 4.11. However,  it is true if the additional condition that $ \mathcal {B} =  \mathcal {B'} $ holds  which is shown  in the following Theorem 4.14.\\
 
 \textbf{Theorem  4.14 :} If $A$ is $s g^*$-closed  and $ s\lambda^*$-closed  and satisfies the condition $ \mathcal {B} =  \mathcal {B'} $, then $A$ is a semi-closed set.

\begin{proof}
Suppose $A$ is $ sg^*$-closed  and $ s\lambda^*$-closed  satisfying the condition $ \mathcal {B} =  \mathcal {B'} $. Since $ A $ is $s g^* $-closed,  there is a semi-closed set $F$ containing $A$ such that $ F\subset  sA_\tau^\wedge $.  Now $ A\subset F $ implies $ \overline{A}_s\subset\overline{F}_s = F\subset s A_\tau^\wedge $.  Again since $A$ is $s \lambda^*$-closed,  $A=sA_\tau^\wedge \cap\overline{A}_s=\overline{A}_s$..............(1). 

 Now $A$ is $s g^*$-closed , so $ \overline{A_s^*}=A$, therefore $ (X - A) \in \mathcal{B}' $.  Since $ \mathcal {B} =  \mathcal {B'} $, $ (X - A) \in \mathcal{B} $ which implies  $\overline{A}_s $ is semi-closed. Therefore from (1)   we have $A$ is a semi-closed set.  
\end{proof}

\textbf{Remark  4.15 :} Union of two   $ s\lambda^*$-closed  sets may not be  $ s\lambda^*$-closed  set as revealed in the under noted Example 4.16.\\

\textbf{Example  4.16 :} Suppose $ X=R-Q $ and $ \tau = \{X, \emptyset, G_i \}$ where $ \{G_i\} $ are the countable subsets of $ X- \{\sqrt{2}\} $. Then $ (X, \tau) $ is a space, but not a topological space. Suppose $ A = X - \{\sqrt{2}, \sqrt{3}, \sqrt{5}\} $ and $ B = X - \{\sqrt{2}, \sqrt{7}, \sqrt{11}\} $. Then $ X- A $ is semi-open set since $ \{\sqrt{3}, \sqrt{5}\} $ is an open set and $ \{\sqrt{3}, \sqrt{5}\}  \subset X - A \subset \overline{\{\sqrt{3}, \sqrt{5}\}} $. Similarly $ X - B $ is also semi-open. Therefore $ A $ and $ B $ are semi closed sets. Now $ sA_\tau^\wedge \cap \overline{A}_s = A $ and $ sB_\tau^\wedge\cap \overline{B}_s = B $. So $ A $ and $ B $ are $ s\lambda^* $-closed sets.  Again take $ C = A \cup B = \{X - \{\sqrt{2}, \sqrt{3}, \sqrt{5}\}\}\cup \{X - \{\sqrt{2}, \sqrt{7}, \sqrt{11}\}\} = X - \{ \sqrt{2}\} $. Since $ \{\sqrt{2}\}$ is not  semi-open, $ X - \{\sqrt{2}\} $ is not semi-closed. So $ \overline{C}_s = X $.  Since the closure of any proper open set $ D $ contained in $ C $ will be $ D \cup\{\sqrt{2}\} $, so $ D \subset C \not\subset \overline{D} $ since $ C $ is not countable where as $ \overline{D} $ is countable, this implies $ X - \{\sqrt{2}\} $ is not semi-open. Therefore $ sC_\tau^\wedge = X $. So $ \overline{C}_s \cap sC_\tau^\wedge = X  \not = C $, hence $ C $ is not $ s \lambda^*$-closed  set.  \\

\textbf{Note 4.17 }: In view of the  Example 4.16, it follows that intersection of two $s \lambda^* $-open sets may not be $s \lambda^* $-open.\\

\textbf{Theorem 4.18 }: A subset $ A $ of $ X $   is  $s \lambda^* $-open if and only if  $ A = N\cup sInt(H) $ where $ N $ is a $s \vee_\tau $-set  and $ H $ is a subset of $X$.
\begin{proof}
Let $A$ be a $s \lambda ^* $-open set. Then $ X - A $ is a $s \lambda^* $-closed set. Then, $ X - A = K\cap\overline{P}_s $, where $K$ is a $s \wedge_\tau $-set and $ P\subset X $. Therefore, $ A = (K\cap\overline{P_s})^c = K^c\cup (\overline{P}_s)^c = (X - K)\cup (X - \overline{P}_s)$, where $ X - K $ is a $s \vee_\tau $-set and $ (X -\overline{P}_s) \subset X $. Take $ N = X - K \subset A $ and $ H = X-P $. So $ P=X-H $ implies $ \overline{P}_s = \overline{(X-H)}_s $ implies $ X-\overline{(X-H)}_s=X-\overline{P}_s $. So by Theorem 4.4, $ sInt(H)=X-\overline{(X-H)}_s=X-\overline{P}_s$. Therefore $ A = N\cup sInt(H) $.

 Conversely,  Suppose  $ A=N\cup sInt(H) ,  N $ is a $s \vee_\tau $-set and $ H\subset X $.  Therefore $ X- A=N^c\cap (sInt(H))^c  =(X-N)\cap (X-sInt(H))$, where $ N^c = K $ (say) is $s \wedge_\tau$-set. Take $ P=X-H $. By  Theorem 4.4, $ X-\overline{(X-H)}_s= sInt(H)$ implies $ \overline{(X-H)}_s=X-sInt(H)$. Therefore $ X - A=K\cap \overline{P}_s $ implies $ X - A $ is a $s \lambda^* $-closed set and hence $A$ is a $s \lambda ^* $-open set.
\end{proof}

\textbf{Theorem 4.19 :} A subset $ A $ of $ X $   is  $s \lambda^* $-open if and only if  $ A = sA_\tau^\vee \cup sInt(A) $.
\begin{proof}
Suppose $ A $ is  $s \lambda^* $-open, Then $ A = N\cup sInt(H) $ where $ N $ is a $s \vee_\tau $-set  and $ H $ is a subset of $X$. Since $ N\subset A $ then $ sN_\tau^\vee\subset sA_\tau^\vee $ and since $ sInt(H) \subset A $, then $ sInt(sInt(H))\subset sInt(A) $. So $ A=N\cup sInt(H)=sN_\tau^\vee\cup sInt(H)\subset s A_\tau^\vee \cup sInt(H)=sA_\tau^\vee\cup sInt(sInt(H))\subset sA_\tau^\vee\cup sInt(A) $. Again $s A_\tau^\vee\subset A $ and $ sInt(A) \subset A $ then $sA_\tau^\vee\cup sInt(A)\subset A $. Therefore we get $ A=sA_\tau^\vee\cup sInt(A) $.

Conversely,  suppose  $ A=sA_\tau^\vee\cup sInt(A) $.  Since $ sA_\tau^\vee $ is a $s \vee_\tau $-set and $ A\subset X $, then $A$ is $s \lambda^* $-open.
\end{proof}

\textbf{Remark 4.20 :} Clearly $s \vee_\tau $-sets are $s \lambda ^* $-open sets and semi-open sets are $s \lambda^* $-open sets. On the other hand if $A$ is $s \lambda^* $-open, $sg^*$-open,  and $ \mathcal {B} = \mathcal {B}'$   then $ A $ is semi-open.\\

In  \cite{MS} it is seen that the collection of all $ \lambda_\mu $-open sets forms a generalised topology $ \mu $ on $ X $, but in view of the Note 4.17, the collecion of $s \lambda^* $-open sets does not form a space structure $ \sigma $ on $ X $.  \\

\section{\bf  \textit Semi-$T_\frac{\omega}{4}$ space, \quad \textit Semi-$T_\frac{3\omega}{8}$ space,    \quad\textit  Semi-$T_\frac{5\omega}{8}$ space  :}

\textbf{Definition 5.1} :  A space $ (X,  \tau) $  is called semi-$T_\frac{\omega}{4}$ if for every finite subset $P$ of $X$ and for every $ y\in X - P $,  there exists a set $A_y$ containing $P$ and disjoint from $\{y\}$  such that $A_y$ is either semi-open or semi-closed.\\

\textbf{Theorem  5.2} :  A space $ (X,  \tau) $ is semi-$ T_\frac{\omega}{4}$ if and only if every finite subset of $X$ is $s \lambda^*$-closed.  

\begin{proof}
Suppose $ (X,  \tau) $ is semi-$ T_\frac{\omega}{4}$ space  and $P$ is a finite subset of $X$.  So for every $ y\in X - P $ there is a set $ A_y $ containing $P$ and disjoint from $ \{y\} $ such that $ A_y $ is either semi-open or semi-closed.  Let $K$ be the intersection of all such semi-open sets $ A_y $ and $F$ be the intersection of all such  semi-closed sets  $A_y$ as $ y $ runs over $ X -P $.  Then $ K=sK_\tau^\wedge$ i.e. $K$ is a $ s\wedge_\tau $-set  and $ \overline{F}_s=F$. Therefore $ K\cap F=K\cap \overline{F}_s = P$. So $ P $ is $ s\lambda^*$-closed. 

Conversely,  let $P$ be a finite set. So by the condition it is $s \lambda^* $-closed. Let  $ y\in X - P $. Then $ P=K\cap \overline{P}_s $, by  Theorem 4.6 (iii) where $K$ is a $s \wedge_\tau $-set. If $ y\not\in \overline{P}_s $ then there exists a semi-closed set $ A_y $ containing $P$ such that $ \{y\} \not\subset A_y $.  Again if $ y\in \overline{P}_s $, then $ y\not\in K $. Then there exists some semi-open set $ A_y $ containing $K$ such that $ y\not\in A_y $ and $ A_y $ contains $ P $ also.  Hence $ (X,  \tau) $ is semi-$ T_\frac{\omega}{4}$ space.
\end{proof}

\textbf{Theorem  5.3} : A space $ (X,  \tau) $ is semi-$ T_0 $ if and only if every singleton of $X$ is $s \lambda^*$-closed. 

\begin{proof}
Let the  space $ (X,  \tau) $ be  semi-$ T_0$ and $ x\in X $.  Take a point $ y \in X - \{x\} $. Since $ (X,  \tau) $ is semi-$ T_0 $, there exists a set $ A_y $ containing $ \{x\} $ such that $ y\not \in A_y $ and $ A_y $ is either semi-open or semi-closed. Let $K$ and $P$ be the intersection of all such semi-open sets $ A_y $ and all such semi-closed sets $ A_y $ respectively as $ y $ runs over $ X-\{x\} $. Then $ \overline{P}_s=P $ and $ K=sK_\tau^\wedge $. Therefore $ \{x\} =K\cap\overline{P}_s$. This implies that $ \{x\}$ is a $s \lambda^*$-closed  set by Definition. 

Conversely,  let  $ \{x\} $ be $ s\lambda^*$-closed for each $ x\in X $. Then $ \{x\} = s\{x\}_\tau^\wedge\cap \overline{\{x\}}_s$  by Theorem  4.6 (iv).  Suppose $ y\in X -\{x\} $. If $ y\not\in \overline{\{x\}}_s $,  then there exists a semi-closed set $P$ containing $ \{x\} $ such that $y\not\in P$. If $ y\in \overline{\{x\}}_s $ then $  y\not \in s\{x\}_\tau^\wedge$. Therefore there exists a semi-open set $ V $ containing $ \{x\} $ such that $ y\not\in V $.   Hence the space $ (X,  \tau) $ is semi-$ T_0 $. \end{proof}

\textbf{Theorem 5.4} :    Every semi-$T_\frac{\omega}{4}$  space  is a semi-$T_0$ space.
\begin{proof}
Suppose the space $ (X,  \tau) $ is semi-$T_\frac{\omega}{4}$ and  $ x , y $ are two distinct points in $ X $.  Since the space is semi-$T_\frac{\omega}{4}$, then for every $ y\in X -\{x\} $ there exists a set $A_y$ containing $ x $ for each $ x $ such that  $ A_y$ is either semi-open or semi-closed and $ y\not\in A_y $. This implies that $ (X,  \tau) $ is semi-$T_0$ space.  
\end{proof}

\textbf{Definition 5.5} : A space$ (X,  \tau) $  is called semi-{$T_\frac{3\omega}{8}$ space if for every countable subset $P$ of $X$ and for every $ y\in X - P$,  there exists  a set $ A_y $ containing $P$ and disjoint from $ \{y\} $ such that $ A_y $ is either semi-open or semi-closed.\\

 \textbf{Theorem  5.6} :  A space$ (X,  \tau) $ is semi-{$T_\frac{3\omega}{8}$ space if and only if every countable subset of $X$ is $ s\lambda^*$-closed.
 Proof is similar to the proof of Theorem 5.2,  so is omitted.\\

 \textbf{Definition 5.7 :} A space$ (X,  \tau) $  is called semi-{$T_\frac{5\omega}{8}$ space if for any subset $ P $ of $ X $ and for every $ y\in X - P$,  there exists  a set $ A_y $ containing $P$ and disjoint from $ \{y\} $ such that $ A_y $ is either semi-open or semi-closed.\\

 \textbf{Theorem  5.8  :} A space$ (X,  \tau) $ is semi-{$T_\frac{5\omega}{8}$ space if and only if for every subset $ E $ of $ X $ is $s \lambda^*$-closed .
  Proof is parallel to the proof of Theorem 5.2,  so is omitted.\\

 Note that {$T_\frac{5\omega}{8}$-space does not imply the condition $ \mathcal {B} = \mathcal {B}' $ where $ \mathcal {B} $  and $ \mathcal {B'} $ are as in Theorem 3.31.\\

 \textbf{Remark  5.9  :} It follows from the Theorems 4.12,  5.8,  5.6,   5.2 that every semi-$ {T_\omega}$  space is semi-{$T_\frac{5\omega}{8}$ space and semi-{$T_\frac{5\omega}{8}$ space is semi-{$T_\frac{3\omega}{8}$ space and semi-{$T_\frac{3\omega}{8}$ space is  semi-$ T_\frac{\omega}{4}$ space. 
 
 However, the converse of each implication may not be true which are substantiated respectively in the undermentioned examples.\\
 
  \textbf{Example 5.10 :} Example of a  semi-{$T_\frac{5\omega}{8}$ space which is not semi-$ T_\omega$ space.
 
  Let $X=R-Q$,  $ \tau = \{X,  \emptyset,  G_i\}$,  where $ \{G_i\}$  are the all countable subsets of $X$. So $ (X,  \tau) $ is a space but not a topological space.  Let $ A $ be a  subset of $ X $. Then $\overline{ A}_s=A $. So $ sA_\tau^\wedge \cap \overline{ A}_s= A $.   This implies that $ A $ is $s \lambda^*$-closed. Therefore any subset of $ X $ is $s \lambda^*$-closed which in turn $ (X, \tau) $ is semi-{$T_\frac{5\omega}{8}$ space. Now let $ B $ be the set of all irrationals in $ (0, \infty) $. Then $ B $ is  $sg^*$-closed but not semi-closed since $ X $ is the only semi-open set containing $ B $ and $ X-B $ is not semi-open. This implies the space is not semi-$ T_\omega$ space.\\

 \textbf{Example 5.11 :} Example of a semi-$T_\frac{3\omega}{8}$ space which is not semi-$ T_\frac{5\omega}{8}$ space.
 
  Suppose  $ X=R-Q $,  and $ \tau=\{X,  \emptyset,  G_i\} $ where $\{ G_i\} $ are all countable subsets of $ X - \{\sqrt{2}\}$. Then $ (X,  \tau) $ is a space but not a topological space.  Take any countable subset $A\subset  X $. Then if  $\sqrt{2}\not\in A$, $A$ is an open set, so $ A $ is semi-open which implies that $ A $ is  $s \lambda^*$-closed. Again if $\sqrt{2}\in A, A $ is not open. But since $ A - \{\sqrt{2}\} $ is an open set  and $ (A - \{\sqrt{2}\})\subset A\subset \overline{A - \{\sqrt{2}\}}= A $, then $ A $ is semi-open. Therefore in both cases $ sA_\tau^\wedge=A $. Hence $ A $ is $s\lambda^*$-closed.    So $ (X,  \tau) $ is a semi-{$T_\frac{3\omega}{8}$ space.
  
   Now suppose $ B = X - \{\sqrt{2}\} $ which is an uncountably infinite subset of $X$. Therefore $ B $ is not an open set. Since the     closure of any proper open set $ D \subset B $ will be $ D\cup \{\sqrt{2}\} $ and $ D\subset B \not \subset\overline{D}$, therefore $ X - \{\sqrt{2}\}$ is not semi-open. So $ sB_\tau^\wedge=X $.  Since $ \{\sqrt{2}\}$ is not semi-open, $ X - \{\sqrt{2}\} $ is not semi-closed, so $ \overline{B}_s = X $. Then  $s B_\tau^\wedge\cap \overline{B}_s = X\not=B $. Therefore   $ B $ is not $s \lambda^*$-closed. So $ (X,  \tau )$ is not semi-$ T_\frac{5\omega}{8}$ space.\\

\textbf{Example  5.12  :} Example of a  semi-$T_\frac{\omega}{4}$ space which is not semi-$ T_\frac{3\omega}{8}$ space.

Suppose  $ X=R-Q , X^{*}=\{2\}\cup X,  \tau=\{\emptyset, X^{*},  \{2\}\cup (X - A): A\subset X, A $ is finite\}.  Therefore  $ (X^*,  \tau) $ is a  topological space so a space also.  Take any finite  subset $E \subset X^{*} $,  we get the following observations:

(i)  if $ 2 \in E, sE_\tau^\wedge=\cap\{ (X - \{\alpha\})\cup\{2\} ,\alpha \in X - E\}=\{2\}\cup E =E $. So $ sE_\tau^\wedge \cap \overline{E}_s=E $ which implies that $ E $ is a $s \lambda^*$-closed set .

(ii) if $ 2 \not \in E ,  E$ is a closed set which implies that $ E $ is a semi-closed set. Hence $E$ is $s\lambda^*$-closed.

So $ (X^{*},  \tau) $ is a semi-$T_\frac{\omega}{4}$ space.

Suppose $Y$ is a countably infinite subset of $X$. So $2\not\in Y $ and  $ Y $ itself is not a semi-open set. Therefore $ sY_\tau^\wedge=\{2\}\cup Y $. $ Y $ is not a closed set, since $ Y $ is not a finite set. Also  $ X^*-Y $ is not semi-open because there is no non-empty open set contained in $ X^*-Y $ whose closure is non-empty, $ Y $ is not itself a semi-closed set. We will prove now that no semi-closed set is there containing $ Y $ except $ X^* $. If a proper semi-closed set $ B $ contains $ Y $, then $ X^*-B $ would be semi-open. But any open set contained in $ X^*-B $ is the null set $ \emptyset $ whose closure is $ \emptyset $ which contradicts the definition of semi-open set. Here closed sets are finite. Therefore $ \overline{Y}_s=X^{*} $. Thus $ \overline{Y}_s\cap sY_\tau^\wedge=[\{2\}\cup Y]\cap X^* =\{2\}\cup Y\not=Y $. Therefore $ Y $ is not $s \lambda^{*} $-closed. Hence the result.\\

\textbf{Definition  5.13}: A space $ (X, \tau) $ is said to be semi-$ R_0  $ if every semi-open set contains the semi-closure of each of its singleton.\\

\textbf{Theorem 5.14} :  If $ (X,  \tau )$ is  semi-$T_0 $  then for every pair of distinct points $ p,q\in X $,  either $ p\not \in\overline{\{q\}}_s $  or $ q\not \in \overline
{\{p\}}_s $.
\begin{proof}
Let $ (X,  \tau) $ be  semi-$T_0 $ and $ p, q\in X, p \not= q $. Since $ X $ is  semi-$T_0 $ space, there exists a semi-open set $ U $ which contains only one of $ p, q $. Suppose that $ p\in U $ and $ q \not \in U $. Then the semi-open set $ U $ has an empty intersection with $ \{q\} $. Hence $ p\not\in \{q\}_s'$. Since $ p \not = q,  p\not\in\overline{\{q\}}_s$. Similarly if $ U $ contains the point $ q $ then $ q\not\in \overline{\{p\}}_s $. 
\end{proof}

\textbf{Theorem  5.15} :  A space $ (X,  \tau) $ is semi-$ T_1 $ space if and only if it is semi-$ T_0 $ and semi-$ R_0 $. 

\begin{proof}
Let  $ (X,  \tau) $ be semi-$ T_1 $ space.  Obviously then it is semi-$  T_0 $.  Suppose  $ A $ is semi-open and  $ x \in A $. Let $ y\in X , y\not= x $. Since the space is semi-$ T_1 $,  there is semi-open set $ V $ such that $ y\in V $  and $ x\not\in V $. Hence $y$ cannot be a semi-limit point of $ \{x\} $, rather no point lying outside $ \{x\} $ can be a semi-limit point of $ \{x\} $,  so $ \{x\}_s'\subset  \{x\} $. Hence $ \overline{\{x\}}_s =\{x\}\subset A $. So $(X,  \tau) $ is semi-$ R_0 $ space.

Conversely, let $ (X,  \tau) $ be semi-$ T_0 $ and semi-$ R_0 $. So for $ x,  y \in X $ and $ x\not=y $,  either $ x\not\in \overline{\{y\}}_s $ or $ y\not\in \overline{\{x\}}_s $, by Theorem 5.14. Suppose $ x\not\in \overline{\{y\}}_s $,  then there exists a semi-closed set $F$ containing $y$ such that $ x \not\in F $. Therefore $ x\in X - F$,  a semi-open set and $ y \not \in X - F $. Since the space is also semi-$ R_0 $, by definition 5.13,  $\overline{\{x\}}_s \subset X - F $. So $ \overline{\{x\}}_s\cap F=\emptyset $ which implies that $ \overline{\{x\}}_s\cap \{y\}=\emptyset $. Hence $ y \not\in \overline{\{x\}}_s $.  As $ x\not =y $ and  $y \not\in \overline{\{x\}}_s, y$ is not a semi-limit point of $ \{x\} $. Therefore there exists a semi-open set $V$ containing $y$, but $ x \not \in V 
$. Since $ x,  y \in X $ and $ x\not=y $,  we get two semi-open sets $ X - F$ and $V$ containing $x , y$ respectively and $ y \not \in X - F $ and $ x \not \in V $.  Thus $ (X,  \tau) $ is semi-$ T_1 $ space.
\end{proof}

\textbf{Theorem 5.16} :  A space $ (X,  \tau) $ is semi-$ T_1 $ if and only if   every  singleton is $s \wedge_\tau $-set.   

\begin{proof}
Let the space $ (X,  \tau) $ be semi-$ T_1$.    By Theorem 5.15, $ (X, \tau) $ is   semi-$  T_0 $ and  semi-$ R_0 $. Since $ (X,\tau)$ is semi-$ T_0 $,  every singleton is $ s\lambda^*$-closed , by Theorem 5.3. Suppose $ x\in X $,  then $ \{x\} $ is $s \lambda^*$-closed , so $ \{x\}=s\{x\}_\tau^\wedge \cap \overline{\{x\}}_s $,  by Theorem 4.6 (iv). We claim that $ \{x\}= s\{x\}_\tau^\wedge $.  If not, there exists $ y \in  s\{x\}_\tau^\wedge - \{x\}$. So $ y\not\in \overline{\{x\}}_s $ and hence there is a closed set $F$,  $ F \supset \{x\}$ such that $ y \not \in F $. Therefore $ y\in X - F $, a semi-open set. Again since $ (X, \tau) $ is  semi-$ R_0,  \overline{\{y\}}_s \subset X - F $. Thus $\overline{\{y\}}_s\cap F=\emptyset $.  Since $ x \in F $, then $ x\not\in\overline{\{y\}}_s$. Therefore there exists a semi-open set $V$ containing $x$ but $ y \not \in V $.  This implies that $ y \not\in s\{x\}_\tau^\wedge $,  a contradiction.  Hence $ \{x\} $ is  a $s \wedge_\tau  $-set.

Conversely, let $ x,y \in X $ and $ x\not=y $.   So $ y \not \in \{x\} $.  By supposition $ \{x\} $ and $ \{y\} $ are $s \wedge_\tau $-sets.  So $ \{x\}=s\{x\}_\tau^\wedge$, therefore $ y \not \in s\{x\}_\tau^\wedge $. So there exists a semi-open set $ V' $ such that $x \in V' $,  but $ y \not \in V' $. Similarly, since $ \{y\}=s\{y\}_\tau^\wedge$ , there exists a semi-open set $ U' $ such that $ y\in U' $ and  $ x \not \in U' $.  Hence $x , y$ are weakly semi-separated by semi-open sets $ V' $ and $ U'$  respectively.
So $ (X,  \tau) $  is semi-$ T_1  $space.  
\end{proof}

\textbf{Corollary  5.17:}  If every subset of $X$ is $s \wedge_\tau $-set then $ (X, \tau) $ is semi-$ T_1 $ space.

The proof follows directly from Theorem 5.16.\\

\textbf{Note:  5.18  :}  The converse of the Corollary 5.17 is not true as revealed in the following Example 5.19. However the converse is true   if   additional conditions are imposed as given in Theorem 5.24. 
Also note that the converse part is true in a $\mu $-space  \cite{MS}.\\

\textbf{ Example 5.19 :}  Let $ X = R-Q $,    $ \{G_i\} $ be the all countable subsets of $X$ and $ \tau $ = $\{X,  \emptyset, G_i\}$. Therefore $ (X, \tau) $ is a space but not a topological space. Since every singleton in $( X, \tau )$ is $ s\wedge_\tau $-set, the space is $ T_1 $ by Theorem 5.16. Let $ B $ be the set of all irrationals in $ (0,  1) $. Then $ B $
is not a $s \wedge_\tau$-set as $ s B_\tau^\wedge=X\not=B $. Hence the result follows.\\

\textbf{Definition  5.20 :}  A space $ (X, \tau) $ is said to be semi-symmetric if $ x,  y \in X,  x \in \overline{\{y\}}_s \Rightarrow y\in \overline{\{x\}}_s $\\ 

\textbf{Definition  5.21 :}  A space $ (X, \tau) $ is said to be strongly semi-symmetric if each singleton in $  X  $  is $s g^*$-closed.\\

\textbf{Theorem  5.22 :} A strongly semi-symmetric space $ (X,  \tau) $ is semi-symmetric.

\begin{proof}
Suppose $x \in \overline{\{y\}}_s$, but $ y\not \in\overline{\{x\}}_s $ for $ x, y \in X $. Then there is a semi-closed set $ F \supset \{x\} $ such that $ y \not \in F $. So $\{y\} \subset F^c $ where $ F^c $ is semi-open. Since $ \{y\} $ is $ sg^* $-closed, there is a semi-closed set $ F' \supset \{y\} $ such that $ F' \subset F^c $. So $ \overline{\{y\}}_s \subset F' $ and hence $ x \in \overline{\{y\}}_s \subset F' \subset F^c $, a contradiction, since $ x \in F $.  Hence the result.
\end{proof}

But the converse may not be true as shown by the undernoted Example 5.23.\\

\textbf{Example  5.23 :} Suppose $ X = R - Q $ and $ \tau = \{X, \emptyset, G_i \}$ where $ \{G_i\}$ are the countable subsets of $ X $. So $(X, \tau)$ is a space but not a topological space. Suppose $ x, y \in X $.   If $ x \not = y $ then $ x\in \overline{\{y\}}_s $ implies that $ x $ is a semi-limit point of $ \{y\} $. But $ \{x\} $ is a semi-open set contains $ x $ which does not intersect $ \{y\} $. So $ x $ cannot be a  semi-limit point of $ \{y\} $. So we must have $ x = y $. So in that case $ y\in \overline{\{x\}}_s $. Hence $ (X, \tau) $ is  semi-symmetric. But for each singleton $ \{x\}$ in $ X ,   s\{x\}_\tau^\wedge = \{x\} $ and $ \{x\} $ is not semi-closed. So no singleton is $ sg^* $-closed. Hence the result.\\

\textbf{ Theorem  5.24  :}  If $ (X,  \tau) $ is a strongly semi-symmetric semi-$ T_1 $ space and satisfies the condition $ \mathcal {B} = \mathcal {B}' $ where $ \mathcal {B} $  and $ \mathcal {B'} $ are as in Theorem 3.31, then every subset of X is a $s \wedge_\tau $-set.
\begin{proof}

Let $ (X,  \tau) $ be a strongly semi-symmetric semi-$ T_1 $ space satisfying the condition $ \mathcal {B} = \mathcal{B}' $ and let $ A\subset X$. Let $ x\in X - A $. Then $ \{x\}$ is $s g^*$-closed. Since $ (X,\tau) $ is semi-$ T_1 $ space,  $ \{x\}  $ is a $s \wedge_\tau $-set by Theorem 5.16, so a $s \lambda^*$-closed set.  Therefore $ \{x\} $ is a semi-closed set by Theorem 4.14.  Therefore $ X - \{x\} $ is a semi-open set containing A.  So $ sA_\tau^\wedge=\cap \{X - \{x\}, x\in X - A\} = A $ which implies $A$ is a  $ s\wedge_\tau $-set.
\end{proof}

\textbf{ Remark 5.25:} If the space  $ (X,  \tau) $ is a strongly semi symmetric semi-$ T_1 $ space and satisfies the condition $ \mathcal{B} = \mathcal {B}' $, then it follows from Theorem 5.24 that union and intersection of two $s \lambda^* $-closed sets are $ s\lambda^* $-closed sets.\\

\textbf{Definition  5.26} : A space $ (X,  \tau) $ is said to be  weak semi-$ R_0$ space  if every $ s\lambda^*$-closed singleton is a $s \wedge_\tau $-set. \\

\textbf{Theorem 5.27 } : Every semi-$ R_0 $  space  is a weak semi-$ R_0 $ space.
\begin{proof}
Let the space be semi-$ R_0 $ and $ x\in X $ be such that $ \{x\} $ be $s \lambda^* $-closed. Then $ \{x\}=s\{x\}_\tau^\wedge\cap\overline{\{x\}}_s$ by  Theorem 4.6 (iv). We claim that $\{x\}  $ is a $s \wedge_\tau $-set i.e.  $ \{x\}=s\{x\}_\tau^\wedge$. If not,  i.e if $  \{x\}\not=s\{x\}_\tau^\wedge  $, then there exists $  y\in s\{x\}_\tau^\wedge -\{x\}$. Therefore $ y\not\in\overline{\{x\}}_s$,  so there is a semi-closed set $F,   F\supset\{x\}  $  such that $ y\not\in F $. This implies that $ y\in X - F $,  a semi-open set.  Since$ (X,\tau) $ is semi-$ R_0  $ space,  $ \overline{\{y\}}_s\subset X - F$. Therefore  $\overline{\{y\}}_s\cap F=\emptyset$.  Since $ x\in F  $,  $ x\not\in\overline{\{y\}}_s=\{y\}\cup \{y\}_s'$ where $ \{y\}_s' $ denotes the set of semi-limit points of $ \{y\} $. Therefore there exists a semi-open set $ V \supset\{x\} $ but $y\not\in V$,  since $x\not=y$ and $ x$ is not also the semi-limit point of $ \{y\} $ . This implies that $ y\not\in s\{{x}\}_\tau^\wedge $,  a contradiction. Hence $ (X,  \tau) $ is weak semi-$ R_0 $ space.
\end {proof}

\textbf{ Theorem 5.28} : For a space $ (X, \tau)$ , the following statements are equivalent:

(1)   $ (X, \tau)$  is semi-$ T_1 $ 

(2)   $ (X, \tau)$  is semi-$ T_0 $ and semi-$ R_0 $

 (3)  $ (X, \tau)$ is semi-$ T_0 $  and weak semi-$ R_0 $.
 
 \begin{proof}
 $(1)\Rightarrow (2)$  follows from Theorem  5.15
 
 $(2) \Rightarrow  (3)$  follows from Theorem  5.27
       
  $(3)   \Rightarrow  (1)$ :  Let $ (X, \tau)$  be semi-$ T_0 $  and weak semi-$ R_0 $ and $ \{x\}\subset X $.  So by Theorem 5.3,  $ \{x\}$  is $s \lambda^*$-closed . Again  $ (X, \tau)$ is weak semi-$ R_0 $,  $\{x\} $ is $s \wedge_\tau$-set. By Theorem 5.16,  $(X, \tau)$ is  semi-$ T_1 $.
  \end{proof}
  
  \textbf{ Theorem 5.29 :} If $ (X ,\tau) $ is a strongly semi-symmetric semi-$ T_1 $ space and satisfies the condition $ \mathcal {B} = \mathcal {B}' $, then it is semi-$ T\omega $ space.
  \begin{proof}
  Let $(X, \tau)$ be strongly semi-symmetric semi-$ T_1 $ space satisfying the condition $ \mathcal{B} = \mathcal {B}'$ and let $ A $ be a $s g^*$-closed  set. Then $\overline{A_s^*}=A $, so $\overline{A_s^*}$ is  $s g^*$-closed. Therefore $ X-A\in \mathcal {B}' $.  Since $\mathcal{B} = \mathcal {B}'$, so $ X-A\in \mathcal{B} $ which implies that $ \overline{A}_s$ is semi-closed. We wish to prove that $ \overline{A}_s = A $.  If not, let $ x\in\overline{A}_s - A $.  Since $ (X, \tau) $ is semi-$ T_1  $ and strongly semi-symmetric space, $ \{x\} $ is  $s g^*$- closed  and $s \lambda^*$-closed and since $ \mathcal{B} = \mathcal {B}'
  , \{x\} $ is a semi-closed set by Theorem 4.14.  But by Theorem 3.13,  $\{x\}\not\subset \overline{A}_s - A$. Therefore $ x\in A $ and so $\overline{A}_s=A$ which implies that $A$  is a semi-closed set and hence $ (X, \tau) $ is semi-$ T_\omega $ space.
  \end{proof}
  
  \textbf{Theorem 5.30  :} If  the space $(X, \tau) $ is strongly semi-symmetric,  weak semi-$ R_0 $,  and satisfies the condition $ \mathcal {B} = \mathcal{B}' $, then the following are equivalent:
  
   (1)  $(X, \tau) $ is semi-$ T_0 $
   
   (2)  $(X, \tau) $ is semi-$ T_1 $
   
   (3)  $(X, \tau) $ is semi-$ T_\omega $
   
   (4)  $(X, \tau) $ is semi-$ T_\frac{5\omega}{8} $
   
   (5)  $(X, \tau) $ is semi-$ T_\frac{3\omega}{8} $
   
   (6)   $(X, \tau) $ is semi-$ T_\frac{\omega}{4} $.
   
   \begin{proof}
   
   (1) $  \Rightarrow(2) $: It follows from Theorem 5.28. 
   
   (2) $ \Rightarrow (3) $: It follows from Theorem 5.29.
   
   (3)  $\Rightarrow $ (4):  It follows from Theorem 4.12 and Theorem 5.8.
   
   (4)  $\Rightarrow $ (5):  It follows from Theorem 5.8 and Theorem 5.6
   
   (5) $ \Rightarrow $(6) :  It follows from Theorem 5.6 and Theorem 5.2
   
   (6)$ \Rightarrow $ (1): It follows from Theorem 5.4.   
  \end{proof}

\textbf{Theorem 5.31 :} A semi-symmetric semi-$ T_0 $ space is semi-$ T_1 $. 
\begin{proof}
Let $ (X,  \tau)$ be semi-symmetric semi-$ T_0 $ space and $ a, b \in X,  a\not= b$. Since  $ (X,  \tau)$ is semi-$ T_0 $, by Theorem 5.14 either $ a\not\in\overline{\{b\}}_s $ or $ b\not\in\overline{\{a\}}_s $ must hold. Let $ a\not\in\overline{\{b\}}_s $ . Then  $ b\not\in\overline{\{a\}}_s $.  For if $ b\in \overline{\{a\}}_s $  then it would imply $ a\in\overline{\{b\}}_s $, since $ (X, \tau) $ is semi-symmetric. But this contradicts that $   a\not\in\overline{\{b\}}_s $. Since $ a\not\in \overline{\{b\}}_s$, there is a semi-closed set $ F $ such that $ b\in F $ and $ a\not\in F $. So $ a\in X-F $, a semi-open set and $ b\not\in X-F $. Again since $ b\not\in\overline{\{a\}}_s $, there is a semi-closed set $ P $ such that $ a\in P $ and $ b\not\in P $.  So $ b\in X-P $, a semi-open set and $ a\not\in X-P $. So $ a,  b $ are weakly semi-separated by semi-open sets $ X-F $ and $ X-P $. Hence $( X,  \tau) $ is semi-$ T_1 $.  
\end{proof}

\textbf{Corollary 5.32 :} If every subset is the intersection of all semi-open sets and all semi-closed sets containing it, then each singleton is either semi-open or semi-closed.
 
\begin{proof}
Let $ x\in X $ but $ \{x\} $ is not semi-open. So $ X-\{x\} $ is not semi-closed. Since $ X-\{x\}\subset X $, so $ X-\{x\}$ is the intersection of all semi-open sets and semi-closed sets containing it. Now we see that $ X $ is the only semi-closed set containing $ X-\{x\} $.  But $ X $ is also a semi-open set containing $ X-\{x\} $, so by the condition $ X-\{x\} $ must be semi-open which in turn $ \{x\} $ is a semi-closed set.
\end{proof}

\textbf{Remark  5.33 :} The equivalence of the conditions stated in the Theorem 5.30 can also be deduced in  terms of another assumption shown in Theorem 5.34 and Theorem 5.35.\\

\textbf{ Theorem 5.34 :} If the space $ (X, \tau) $ is semi-symmetric and satisfies the condition $ \mathcal {B} = \mathcal{B}' $, then the following statements are equivalent:

(1)  $(X, \tau) $ is semi-$ T_0 $
   
(2)  $(X, \tau) $ is semi-$ T_1 $
   
(3)  $(X, \tau) $ is semi-$ T_\omega $
   
(4)  $(X, \tau) $ is semi-$ T_\frac{5\omega}{8} $
   
(5)  $(X, \tau) $ is semi-$ T_\frac{3\omega}{8} $
   
(6)   $(X, \tau) $ is semi-$ T_\frac{\omega}{4} $.
\begin{proof}
(1)  $\Rightarrow$ (2): It follows from Theorem 5.31.

(2)  $\Rightarrow$ (3): It follows from Theorem 5.16, Corollary 5.32 and  Theorem 3.31. 

(3)  $\Rightarrow$ (4): It follows from Theorem 4.12 and Theorem 5.8.

(4)  $\Rightarrow$ (5): It follows from Theorem 5.8 and Theorem 5.6.

(5)  $\Rightarrow$ (6): It follows from Theorem 5.6 and Theorem 5.2.

(6)$\Rightarrow$ (1): It follows from Theorem 5.4.
\end{proof}

\textbf{Theorem 5.35 :} If the space $ (X, \tau) $ is semi-symmetric and  the condition (P) given in Theorem 3.29 is satisfied, then the following statements are equivalent:

(1)  $(X, \tau) $ is semi-$ T_0 $
   
(2)  $(X, \tau) $ is semi-$ T_1 $
   
(3)  $(X, \tau) $ is semi-$ T_\omega $
   
(4)  $(X, \tau) $ is semi-$ T_\frac{5\omega}{8} $
   
(5)  $(X, \tau) $ is semi-$ T_\frac{3\omega}{8} $
   
(6)   $(X, \tau) $ is semi-$ T_\frac{\omega}{4} $.
\begin{proof}
(1)  $\Rightarrow$ (2): It follows from Theorem 5.31.

(2)  $\Rightarrow$ (3): It follows from Theorem 3.29. 

(3)  $\Rightarrow$ (4): It follows from Theorem 4.12 and Theorem 5.8.

(4)  $\Rightarrow$ (5): It follows from Theorem 5.8 and Theorem 5.6.

(5)  $\Rightarrow$ (6): It follows from Theorem 5.6 and Theorem 5.2.

(6)$\Rightarrow$ (1): It follows from Theorem 5.4.
\end{proof}

\end{document}